\documentclass[a4paper,12pt,british]{amsart}

\usepackage[latin1]{inputenc}
\usepackage[T1]{fontenc}
\usepackage{amssymb,lmodern,calrsfs}
\usepackage{babel,hyperref}

\newtheorem{theo}{Theorem}
\newtheorem{coro}[theo]{Corollary}
\newtheorem{prop}[theo]{Proposition}
\newtheorem{lemm}[theo]{Lemma}
\theoremstyle{definition}
\newtheorem{rema}[theo]{Remark}
\newtheorem*{rem}{Remark}
\newtheorem{defi}[theo]{Definition}

\newcommand{\cD}{\mathcal{D}}
\newcommand{\cH}{\mathcal{H}}
\newcommand{\setR}{\mathbb{R}}
\newcommand{\Rr}{\overset{\circ}{R}}
\newcommand{\fn}{\mathfrak{n}}
\newcommand{\mf}{\mathfrak}
\newcommand{\mbb}{\mathbb}
\newcommand{\mc}{\mathcal}
\newcommand{\K}{\mc{K}}
\newcommand{\eL}{\mc{L}}

\newcommand{\be}{\begin{eqnarray*}}
\newcommand{\ee}{\end{eqnarray*}}
\newcommand{\beqa}{\begin{eqnarray}}
\newcommand{\eeqa}{\end{eqnarray}}

\DeclareMathOperator{\Ric}{Ric}
\DeclareMathOperator{\Sym}{Sym}
\DeclareMathOperator{\Tr}{Tr}
\DeclareMathOperator{\im}{im}
\DeclareMathOperator{\End}{End}

\author{S. Armstrong}
\address{St Cross College, Oxford OX1 3LZ}
\author{O. Biquard}
\address{Université Paris 6, UMR 7586, Institut de Mathématiques de Jussieu}
\thanks{The second author is supported by the ANR grant 06-BLAN60154-01.}
\title{Einstein metrics with anisotropic boundary behaviour}

\begin{document}

\maketitle

In recent years the relation between complete, infinite volume,
Einstein metrics and the geometry of their boundary at infinity has
been intensively studied, especially since the advent of the physical
AdS/CFT correspondence.

In all the previous examples of this correspondence, the Einstein
metrics at infinity are supposed to be asymptotic to some fixed
model---a symmetric space of noncompact type $G/K$. Here we shall
restrict to the rank one case, where the examples are asymptotically
real, complex or quaternionic hyperbolic metrics. The corresponding
geometries at infinity (``parabolic geometries'' modelled on $G/P$,
where $P$ is a minimal parabolic subgroup of $G$) are conformal
metrics, CR structures or quaternionic-contact structures. In this
article, we introduce a new class of examples, which are no more
asymptotic to a symmetric space. Actually the model at infinity is
still given by a homogeneous Einstein space, which may vary from point
to point on the boundary at infinity.

This phenomenon cannot occur in the most classical examples (real or
complex hyperbolic spaces), because the algebraic structure at
infinity (abelian group or Heisenberg group) has no deformation. But
such deformations exist for the quaternionic Heisenberg group (except
in dimension 7), and even in the 15-dimensional octonionic case. So
these are the two cases on which this article shall focus. In the
parabolic geometry language, these are the two cases where non regular
examples exist.

More concretely, the basic quaternionic example is the sphere
$S^{4m-1}$, with its $(4m-4)$-dimensional distribution $\cD$, and the
octonionic example is the sphere $S^{15}$ with a $8$-dimensional
distribution $\cD$. At each point $x$ of the sphere, there is an
induced nilpotent Lie algebra structure on $\fn_x=\cD_x\oplus T_xS/\cD_x$,
given by the projection on $T_xS/\cD_x$ of the bracket of two vector
fields $X,Y\in \cD_x$. It was proved in \cite{Biq00} that small deformations of
$\cD$, such that $\fn_x$ remains the quaternionic Heisenberg algebra
for all $x$, are boundaries at infinity of complete Einstein metrics
on the ball. This regularity assumption (that is, keeping the
isomorphism type of the algebra $\fn_x$ fixed) is a strong
differential system on $\cD$: it was shown in \cite{Biq00} that such
quaternionic-contact structures exist in abundance, but there is no
octonionic example \cite{Yam93}.

In this article, we relax the regularity assumption in these two
cases. There is a beautiful family of examples, already known in the
literature: the homogeneous Einstein metrics of Heber \cite{Heb98}. In
the upper space model, each hyperbolic space is identified with the
solvable group $S=AN$, with boundary at infinity the Heisenberg group
$N$ (where $G=KAN$ is the Iwasawa decomposition). Then Heber proved
that every deformation of $S$ carries a unique homogeneous Einstein
metric. In particular, we can associate to a deformation of the
nilpotent Lie algebra $\fn$ the homogeneous Einstein metric on the
corresponding solvable group $S=AN$.

\begin{theo}\label{th:main}
  Let $n=4m-1\geq 11$ in the quaternionic case, or $n=15$ in the
  octonionic case. Any small deformation of the $(4m-4)$-dimensional (in
  the quaternionic case) or $8$-dimensional (in the octonionic case)
  distribution of $S^n$ is the boundary at infinity of a complete
  Einstein metric on the ball $B^{n+1}$.

  At each point $x\in S^n$, the Einstein metric is asymptotic to Heber's
  homogeneous metric on the solvable group associated to the nilpotent
  algebra $\fn_x$.
\end{theo}

The meaning of the theorem is that all deformations of the
distributions on the boundary of the rank one symmetric spaces can be
interpreted as boundaries at infinity of Einstein metrics, but maybe
with an anisotropic behaviour (the asymptotics depends on the
direction). This gives new examples in the quaternionic case, for
dimension at least 11, and in the octonionic case.

The asymptotic condition means that there is some inhomogeneous
rescaling of the metric near the boundary point $x$ which converges to
Heber's metric, see remark \ref{rema:asympt-solv} for details.

The relation between the regular examples and the new examples is
perhaps best understood by remembering that the Einstein metrics
associated to quaternionic-contact structures in dimension at least 11
are actually quaternionic-Kähler \cite{Biq02}, so they keep the holonomy
$Sp_mSp_1$ of the hyperbolic space. This condition distinguishes
exactly the regular case:

\begin{coro}
  In the quaternionic case, for $m\geq 3$, the Einstein metric
  constructed by the previous theorem is quaternionic-Kähler if and
  only if the distribution on $S^{4m-1}$ is regular (that is, is a
  quaternionic-contact structure).
\end{coro}

The corollary follows from the fact that the boundary at infinity of a
quaternionic-Kähler metric must be a quaternionic-contact structure
\cite{Biq00}.

There is a similar, but obvious, story in the octonionic case. The
Cayley plane has holonomy $Spin_9$. If the Einstein metric keeps the
$Spin_9$ condition, it is well-known that it is the hyperbolic metric
($Spin_9$ metrics are locally symmetric). On the other hand, a regular
distribution of dimension 8 on $S^{15}$ must be standard. So we have a
(trivial) example of the equivalence of the holonomy condition on the
Einstein metric with the regularity condition on the boundary.

The article has two parts. The first part is algebraic, and consists
in the construction of an approximate Einstein metric near the
boundary at infinity. The new point here is that the model is not
explicit: it is the solution of algebraic equations giving conformal
structures on the distribution $\cD$ and on the quotient
$TS/\cD$. These equations have a nice interpretation in terms of a
stronger geometric structure, a quaternionic (or octonionic) structure
on $\cD$, on which we add a gauge condition which enables to find a
unique solution. This additional structure should be useful in future
applications, in particular if one wishes to work out a
Fefferman-Graham type development of the Einstein metric.

The second part is analytic, and consists in deforming an approximate
Einstein metric into a solution of the equations. This relies
basically on a deformation argument, which requires to understand the
analytic properties of the deformation operator. If one has a good
understanding of the analysis for the models (Einstein metrics on
solvable groups), then one can probably use microlocal analysis to
glue together the inverses of the deformation operator into the
required parametrix. However here we prefer to avoid the analysis on
these solvable groups, since more direct methods give the required
result. Nevertheless, it is clear that the more sophisticated
microlocal analysis may be required in further developments of the
theory.

\section{Algebraic considerations}
Let $V_1$ and $V_2$ be vector spaces of dimensions $4m-4$ and $3$ (in
the quaternionic case) or of dimensions $8$ and $7$ (in the octonionic
case). A formal Levi bracket is an element $\ell$ of $W = \land^2 V_1^*
\otimes V_2$. This bracket makes
\begin{equation}
\mf{n} = V_1 \oplus V_2\label{eq:15}
\end{equation}
into a two-graded nilpotent Lie algebra (as the Jacobi identity is
trivially satisfied). The corresponding Lie group $N$ will then carry
an invariant distribution $\mc{D}$ of same rank as $V_1$. Consider the
Lie bracket $[,]$ on sections of $\mc{D}$. This is a differential
bracket, but the differential part of it only maps into
$\mc{D}$. Hence the map \be
\mc{L} : \land^2 \Gamma(\mc{D}) &\to & \Gamma(TX/\mc{D}) \\
\mc{L} (X,Y) &=& [X,Y]/\mc{D} \ee is an algebraic map, i.e. a section
of $\land^2 \mc{D}^* \otimes (TN/\mc{D})$. If we designate by $L$ the group
$GL(V_1) \oplus GL(V_2)$, then $\mc{D}$ and $TN/\mc{D}$ are bundles
associated to an $L$-principle bundle $E \to N$.

In this set-up, $\mc{L}$ corresponds to an $L$-equivariant map
$f_{\mc{L}}$ from $E$ to $W$. Designate by $E_p$ the fibre of $E$ at
$p \in M$. By construction, $\ell$ is in the image $f_{\mc{L}}(E_p)$ for
all $p$, and this image consists precisely of the $L$ orbit of $\ell$ in
$W$.

Under the identifications $V_1 = \mbb{H}^{m-1}$ and $V_2 =\im
\mbb{H}$, the quaternionic standard Levi bracket $\kappa$ is given by the
choice of a hermitian metric $h$ on $V_1$; in this case, $\kappa$ is simply
the imaginary part of $h$.

Similarly, the standard octonionic bracket (also designated $\kappa$) is
also defined by identifications $V_1 = \mbb{O}$, $V_2 = \im \mbb{O}$
and a choice of $h$.

In general, $\kappa$ is only defined up to $L$-action; but as $\ell$ is
only defined up to $L$-action, we will assume our choice of $\kappa$ is
fixed.

We identify $G_0$ as the stabiliser of $\kappa$ in $L$; this can be seen as
the group that stabilises the quaternionic or octonionic structure. In
the quaternionic case, \be G_0 = \mbb{R}_+^* Sp(1)Sp(m-1) ,
\ee while in the octonionic case, \be G_0 = \mbb{R}_+^* Spin(7). \ee

In general, we consider a manifold $X^n$ of dimension $n=4m-1$ in the
quaternionic case, or $n=15$ in the octonionic case, and a
distribution $\cD\subset TX$ of dimension equal to the rank of $V_1$. At each
point $x$ of $X$, the image in $TX/\cD$ of the bracket $[X_1,X_2]$ of
two vector fields in $\cD$ is an algebraic map, $\mc{L}_x\in \Lambda^2\cD_x^*\otimes(T_xX/\cD_x)$.

The main result of this section is:

\begin{prop}\label{prop:constr-ae}
  There exists a $L$-invariant open set $U\subset W$ (that is an open set of
  $L$-orbits), containing $\kappa$, with the following property. If $X^n$
  has a distribution $\cD$ such that for every $x\in X$ the induced
  bracket $\mc{L}_x\in U$, then there exist metrics $\eta^2$ and $\gamma$ on $TX/\cD$
  and $\cD$, such that, choosing any splitting $TX=\cD\oplus V$, the metric
  \begin{equation}
   g = \frac{dt^2+\eta^2}{t^2} + \frac \gamma t \label{eq:14}
  \end{equation}
  on $\setR_+^*\times X$ is asymptotically Einstein when $t\to0$:
  \begin{equation}
  \Ric(g)=\lambda g + O(t^{\frac 12}) ,\label{eq:16}
  \end{equation}
  where $\lambda=-m-2$ in the quaternionic case, $\lambda=-9$ in the octonionic
  case.  Moreover, this choice of $\eta$ and $\gamma$ is unique, up to the
  conformal transformation: \be (\eta^2, \gamma) \to ( f^2 \eta^2, f \gamma) \ee for $f$
  a strictly positive function $X \to \mbb{R}$.
\end{prop}

In this statement, it is important to note that the asymptotic
behaviour (\ref{eq:16}) does not depend on the choice of splitting
$TX=\cD\oplus V$.

\begin{rema}\label{rema:solvable}
  In the special case where $X$ is the nilpotent group $N$ associated
  to an algebraic bracket $\ell\in \Lambda^2V_1^*\otimes V_2$, and the distribution
  $\cD$ is the associated distribution, then the splitting
  (\ref{eq:15}) gives a canonical choice for $V$. Then the metric
  (\ref{eq:14}) is an invariant metric on the solvable group
  $S=\setR_+^*\ltimes N$. Being asymptotically Einstein when $t\to0$ and
  invariant implies that it is exactly Einstein. This is the metric
  constructed by Heber \cite{Heb98} on $S$. The proof of the
  proposition will give another construction of this metric, at least
  for small deformations of the distribution. Conversely, it also
  follows from our proof that the open set $U$ can be taken equal to
  the set of brackets $\ell$ such that an Einstein metric exists on the
  associated solvable group.
\end{rema}

\begin{rema}\label{rema:asympt-solv}
  In general, at each point $x\in X$ is associated a nilpotent group
  $N_x$ and a solvable group $S_x=\setR_+^*\ltimes N$. We are going to see
  the relation between the asymptotically Einstein metric
  (\ref{eq:16}) and $S_x$. To simplify notation, let us consider only
  the quaternionic case (but the octonionic case is similar). Near $x$
  we choose coordinates $(x_1,\dots,x_n)$ on $X$ such that $\cD_x$ is
  generated by the vector fields $(\frac \partial{\partial x_4},\dots,\frac \partial{\partial
    x_n})$. The distribution $\cD$ is given by the kernel of three
  1-forms, $\eta_1$, $\eta_2$ and $\eta_3$, and we can suppose that at the
  point $x$ one has $\eta_i=dx_i$. Then we consider the homothety
  \begin{equation}
    \label{eq:17}
    h_r(t,x_1,\dots,x_n) = (rt,rx_1,rx_2,rx_3,\sqrt r x_4,\dots,\sqrt r x_n) .
  \end{equation}
 Note $\eta_i=\eta_i^jdx_j$, with $\eta_i^j(0)=\delta_i^j$. Then one has
$$ \bar \eta_i:=\lim_{r\to0}\frac 1{\sqrt r} h_r^*\eta_i = dx_i + \sum_{j,k=4}^n x_k
\frac{\partial\eta_i^j}{\partial x_k}(0) dx_j . $$ The three forms $\bar \eta_1$, $\bar
\eta_2$ and $\bar \eta_3$ are homogeneous, and define exactly the horizontal
distribution of the nilpotent group $N_x$. Denote $\bar \gamma:=\gamma(0)$,
then, when $r\to0$, one obtains the limit
$$ h_r^* g \longrightarrow \frac{dt^2+\bar \eta^2}{t^2} + \frac{\bar \gamma}t . $$
This is an invariant metric on the solvable group $S_x$, and more
precisely it is the homogeneous Einstein metric on $S_x$ mentioned in
the previous remark. This justifies the statement in Theorem
\ref{th:main} that at each point the constructed metric is asymptotic
to the corresponding Heber's metric.
\end{rema}

\begin{proof}[Proof of proposition \ref{prop:constr-ae}]
  The uniqueness comes from remark \ref{rema:asympt-solv} and the
  uniqueness of the homogeneous Einstein metric on $S_x$ proved by
  Heber. Later in the paper, a weaker uniqueness will also be proved.

  We will calculate the Ricci tensor of the metric (\ref{eq:14}) as a
  function of $\gamma$ and $\eta$. The calculation is local, so we can
  choose orthonormal frames $\{\check{X}_i\}$ and $\{\check{Y}_i\}$ of
  $(TX/ \mc{D})$ and $\mc{D}$, respectively.

  On $M = \mbb{R}_+^* \times X$, we can define an orthonormal frame via: \be
  X_0 &=& t \frac{\partial}{\partial t}, \\
  X_i &=& t \check{X}_i, \\
  Y_i &=& \sqrt{t} \check{Y}_i.  \ee Let $O(a)$ denote sections of
  $TM$ whose norm (under $g$) tends to zero at least as last as
  $t^a$. Then we may calculate the Lie brackets of the above frame
  elements: \be
  [X_0, X_i] &=& X_i, \\[0.cm]
  [X_0, Y_i] &=& \frac{1}{2} Y_i, \\[0.cm]
  [X_i, X_j] &=& O(1), \\[0.cm]
  [X_i, Y_j] &=& O(1/2), \\[0.cm]
  [Y_i, Y_j] &=& \eL_{ij} + O(1/2), \\[0.cm]
  \ee where $\eL_{ij} = \eL (Y_i, Y_j)$. In future, we will denote by
  $\eL_{ij}^k$ the $k$ component of $\eL_{ij}$ -- i.e.  \be \eL_{ij}^k
  = g(\eL_{ij}, X_k), \ee and $\eL^k_i$ will be the section of
  $\mc{D}$ defined by: \be g (\eL^k_i, Y_j) = \eL^{k}_{ij}.  \ee

  Now let $\nabla$ be the Levi-Civita connection of $g$. We can
  calculate $\nabla$ by using the Koszul formula: \be
  2g(\nabla_X Y, Z) &=& X \cdot g(Y,Z) + Y \cdot g(X,Z) - Z \cdot g(X,Y) \\
  && g([X,Y],Z) - g([X,Z],Y) - g([Y,Z],X).  \ee Since our frame
  elements are orthonormal, the formula reduces to \be 2g(\nabla_X Y,
  Z) &=& g([X,Y],Z) - g([X,Z],Y) - g([Y,Z],X), \ee giving: \be
  \nabla_{X_0} X_0 &=& 0, \\
  \nabla_{X_0} X_i = \nabla_{X_0} Y_i &=& 0, \\
  \nabla_{X_i} X_0 &=& - X_i, \\
  \nabla_{Y_i} X_0 &=& - \frac{1}{2} Y_i, \\
  \nabla_{X_i} X_j &=& \delta_{ij} X_0 + O(1/2), \\
  \nabla_{X_i} Y_j = \nabla_{Y_j} X_i &=&  -\frac{1}{2} \mc{L}^i_j + O(1/2), \\
  \nabla_{Y_i} Y_j = \nabla_{Y_j} Y_i &=& \frac{1}{2} \eL_{ij} +
  \frac{1}{2} \delta_{ij} X_0 + O(1/2).  \ee

  So in this frame, $\nabla = d + A + O(1/2)$ where $A \in \Gamma(T^*X
  \otimes \End TX)$ is independent of $t$. In detail: \be
  A(X_0) &=& 0, \\
  A(X_i) &:& \left\{
\begin{array}{ccl}
X_0 &\to& - X_i \\
X_j &\to&  \delta_{ij} X_0 \\
Y_j &\to& - \frac{1}{2} \eL^i_j,
\end{array}
\right.
\\ \\
A(Y_i) &:& \left\{
\begin{array}{ccl}
X_0 &\to& - \frac{1}{2} Y_i \\
X_j &\to& - \frac{1}{2} \eL^j_i \\
Y_j &\to& \frac{1}{2} \eL_{ij} + \frac{1}{2} \delta_{ij} X_0.
\end{array}
\right.  \ee In this frame, define $dA(X,Y) = X \cdot A(Y) - Y \cdot
A(X) - A([X,Y])$. Note that differentiating $A$ in the $X_0$ direction
is zero, while differentiating $A$ in the direction of $X_i$ or $Y_i$
picks up a $t$ or $\sqrt{t}$ term, and hence become $O(1/2)$. Thus $d
A(X,Y) = -A([X,Y]) + O(1/2)$.

The curvature $R$ of $\nabla$ is $d A + [A,A]$, which immediately implies that
\be
R_{X_0, X_i} &=& -A(X_i) + O(1/2), \\
R_{X_0, Y_i} &=& -\frac{1}{2} A(Y_i) + O(1/2), \\
R_{X_i, X_j} &=& [A(X_i), A(X_j) ] + O(1/2), \\
R_{X_i, Y_j} &=& [A(X_i), A(Y_j) ] + O(1/2), \\
R_{Y_i, Y_j} &=& - A(\eL_{ij}) + [A(X_i), A(X_j) ] + O(1/2). \\
\ee
The commutator terms are given by:
\be
[A(X_i), A(X_j) ] &:& \left\{
\begin{array}{ccl}
X_0 & \to & 0 \\
X_k & \to & \delta_{ik} X_j - \delta_{jk} X_i \\
Y_k & \to & \frac{1}{4} \left( \eL_{\eL^j_k}^i - \eL_{\eL^i_k}^j \right),
\end{array}
\right.
\\ \\[0.cm]
[A(X_i), A(Y_j) ] &:& \left\{
\begin{array}{ccl}
X_0 & \to & - \frac{1}{4} \eL^i_j \\
X_k & \to & \frac{1}{4} \eL_{\eL_j^k}^i + \frac{1}{2} \delta_{ik} Y_j \\
Y_k & \to &  - \frac{1}{2} \delta_{jk} X_i + \frac{1}{4} \eL_{j\eL^i_k} , \\
\end{array}
\right.
\\ \\[0.cm]
[A(Y_i), A(Y_j) ] &:& \left\{
\begin{array}{ccl}
X_0 & \to & \frac{1}{2} \eL_{ji} \\
X_k & \to & \frac{1}{4} \left( \eL_{j\eL^k_i} - \eL_{i\eL^k_j} \right) + \frac{1}{2} \eL_{ji}^k X_0 \\
Y_k & \to & \frac{1}{4} \left( \eL^{\eL_{ik}}_j - \eL^{\eL_{jk}}_i + \delta_{ik} Y_j - \delta_{jk} Y_i \right).
\end{array}
\right.
\ee

Now we need to take the Ricci-trace of this expression: \be
\Ric_{X_0, X_0} &=& \sum_i g(X_i, R_{X_i,X_0} X_0 ) + \sum_i g (Y_i, R_{Y_i,X_0} X_0 ) \\
&=& \sum_i g(X_i,-X_i) + \sum_i g(Y_i,-\frac{1}{4} Y_i) + O(1/2)\\
&=& \lambda + O(1/2).  \ee Here $\lambda$ is equal to $-3 - (4m-4)/4 =
-m-2$ in the quaternionic case, and $-7 -8/4 = -9 $ in the octonionic
case. The cross-terms of the Ricci curvature all vanish: \be
\Ric_{X_0, X_i} &=& \sum_j g(X_j, R_{X_j,X_0} X_i ) + \sum_i g (Y_j, R_{Y_j,X_0} X_i ) \\
&=& \sum_j g(X_j, \delta_{ij} X_0) + \sum_j \frac{1}{4}\eL_{jj} + O(1/2) \\
&=& O(1/2), \\
\Ric_{X_0, Y_i} &=& \sum_j g(X_j, R_{X_j,X_0} Y_i ) + \sum_i g (Y_j, R_{Y_j,X_0} Y_i ) \\
&=& O(1/2), \\
\Ric_{X_i, Y_j} &=& g(X_0, R_{X_0, X_i} Y_j) + \sum_j g(X_j, R_{X_j,X_i} Y_j ) + \sum_i g (Y_j, R_{Y_j,X_i} Y_j ) \\
&=& O(1/2), \ee the last two expressions vanishing because they are
sums of terms of type $g(X,Y)$ with $X \perp Y$. Next, the $\mc{D} \times
\mc{D}$ term is: \be
\Ric_{X_i, X_j} &=& g(X_0, R_{X_0,X_i} X_j) + \sum_k g(X_k, R_{X_k,X_i} X_j) + \sum_k g(Y_k, R_{Y_k,X_i} X_j) \\
&=& - \delta_{ij} + \delta_{ij} - \sum_k  \delta_{ij} g(X_k, X_k) - \sum_k \frac{1}{2} \delta_{ij}  g(Y_k, Y_k) \\
&& - \sum_k \frac{1}{4} g(Y_k, \eL^i_{\eL^j_k}) + O(1/2).  \ee In the
quaternionic case, this is \beqa \label{qua:one} \Ric_{X_i, X_j} =
\lambda \delta_{ij} + (1-m) \delta_{ij} + \sum_{k=1}^{4m-4}
\frac{1}{4} \eL^i_{k \eL^j_k} + O(1/2).  \eeqa In the octonionic case,
this is \beqa \label{oct:one} \Ric_{X_i, X_j} = \lambda \delta_{ij} -2
\delta_{ij} + \sum_{k=1}^8 \frac{1}{4} \eL^i_{k \eL^j_k} + O(1/2).
\eeqa Finally the $(TX/\mc{D}) \times (TX/\mc{D})$ term is \be
\Ric_{Y_i, Y_j} &=& g(X_0, R_{X_0,Y_i} Y_j) + \sum_k g(X_k, R_{X_k,Y_i} Y_j) + \sum_k g(Y_k, R_{Y_k,Y_i} Y_j) \\
&=& - \frac{1}{4} \delta_{ij} + \sum_k \left(\frac{1}{4} \eL_{i\eL_{j}^k}^k- \frac{1}{2}\delta_{ij} \right) + \frac{1}{4} \left( \delta_{ij} + \sum_k 3 \eL_{ik}^{\eL_{kj}} -\delta_{ij} \right) \\
&& + O(1/2) \\
&=& -\frac{1}{2} \sum_k \delta_{ij} + \frac{1}{4} \sum_k (2
\eL_{ik}^{\eL_{kj}} - \delta_{ij}) + O(1/2).  \ee since \be \sum_k
\eL_{ik}^{\eL_{kj}} = \sum_{kp} \eL_{ik}^p \eL_{kj}^p = \sum_{pk}
\eL_{ip}^k \eL_{pj}^k = \sum_k -\eL_{i \eL_{j}^k}^k.  \ee In the
quaternionic case, the curvature is \beqa \label{qua:two} \Ric_{Y_i,
  Y_j} = \lambda \delta_{ij} + \frac{3}{2} \delta_{ij} + \frac{1}{2}
\sum_{k = 1}^{4m-4} \eL_{ik}^{\eL_{kj}} + O(1/2).  \eeqa In the
octonionic case, it is: \beqa \label{oct:two} \Ric_{Y_i, Y_j} =
\lambda \delta_{ij} + \frac{7}{2} \delta_{ij} + \frac{1}{2} \sum_{k =
  1}^{8} \eL_{ik}^{\eL_{kj}} + O(1/2).  \eeqa
Now $(M,g)$ is asymptotically Einstein if $\Ric_{X_i, X_j} = \lambda \delta_{ij}
+ O(1/2)$ and $\Ric_{Y_i, Y_j} = \lambda \delta_{ij} + O(1/2)$. From now on, we
will use the Einstein summation convention, where any repeated index
is summed over. Then the equations (\ref{qua:one}) and (\ref{qua:two})
imply that in the quaternionic case, we must have:
\beqa \label{qua:con}
\begin{array}{rcl}
\eL_{ij}^k \eL_{op}^q \gamma^{io} \gamma^{jp} &=& 4(m-1) \eta^{kq}, \\
\eL_{ij}^k \eL_{op}^q \gamma^{io} \eta_{kq} &=& 3 \gamma_{jp},
\end{array}
\eeqa
while equations (\ref{oct:one}) and (\ref{oct:two}) imply that in the
octonionic case, we must have: 
\beqa \label{oct:con}
\begin{array}{rcl}
\eL_{ij}^k \eL_{op}^q \gamma^{io} \gamma^{jp} &=& 8 \eta^{kq}, \\
\eL_{ij}^k \eL_{op}^q \gamma^{io} \eta_{kq} &=& 7 \gamma_{jp}.
\end{array}
\eeqa
For an $\ell$ sufficiently close to $\kappa$, these equations can be solved
(see Theorem \ref{infi:ein:theo}), and the solution is unique up to
conformal transformations.
\end{proof}

A natural question is whether, as in the quaternionic-contact case,
the conformal class $(\eta,\gamma)$ comes with a quaternionic structure on
$\cD$ and $TX/\cD$. The same applies for the octonionic-contact
structures, of course. We propose here a construction, where instead
of looking only for a conformal class, one constructs directly a
quaternionic or octonionic structure. As a byproduct, the system
(\ref{qua:con}) or (\ref{oct:con}) is interpreted in a natural way,
see (\ref{al:be:eq}), and existence of a solution is provided.

The automorphism group of these structures is $G_0$, which is
contained in the conformal automorphism group \be G' = \mbb{R}_+^* \times
SO(\eta) \times SO(\gamma) \ee of $(\eta,\gamma)$. Thus it seems that to get the
quaternionic/octonionic structures on the manifold, we need to impose
extra equations beyond (\ref{qua:con}) and (\ref{oct:con}).

These can best be understood by looking at the normality $\partial^*$
operator described in \cite{CartEquiv}, \cite{TCPG} and
\cite{capslo}. It is an algebraic Lie algebra co-differential, which
extends naturally to a bundle operator on associated bundles. If $X$
is a quaternionic- or octonionic-contact manifold, $\K$ the
corresponding Levi-bracket, and $\mc{M}$ is any section of $\land^2
\mc{D}^* \otimes (TX/\mc{D})$. Then $\partial^* \mc{M} = \alpha \mc{M} \oplus \beta \mc{M}$
where \beqa \label{al:be:k}
\begin{array}{rcl}
(\alpha \mc{M})^r_q &=& (\gamma^{jr} \gamma^{ip} \eta_{ko}) (\mc{M}_{ij}^k \mc{K}_{pq}^o), \\
(\beta \mc{M})_r^k &=& - \frac{1}{2}   (\eta_{or} \gamma^{ip} \gamma^{jq}) (\mc{M}_{ij}^k \mc{K}_{pq}^o),
\end{array}
\eeqa
Einstein summation over repeated indexes being assumed. Note that these expressions are invariant under conformal transformations $(\eta, \gamma) \to (f^2 \eta, f \gamma)$.
If we apply $\alpha$ and $\beta$ to $\K$ itself, we get:
\begin{lemm} \label{hom:model}
In the quaternionic-contact case:
\be
(\alpha \mc{K}) &=& 3 Id_{\mc{D}}\\
-2(\beta \mc{K}) &=& 4(m-1) Id_{TX/\mc{D}},
\ee
while in the octonionic-contact case:
\be
(\alpha \mc{K}) &=& 7 Id_{\mc{D}}\\
-2(\beta \mc{K}) &=& 8 Id_{TX/\mc{D}},
\ee
the same numbers as in equations (\ref{qua:con}) and (\ref{oct:con}).
\end{lemm}
\begin{proof}
  Fix $\eta$ and $\gamma$, and pick local orthonormal sections $\{I_1, \cdots, I_p\}$
  of $TX / \mc{D}$, where $p = 3$ in the quaternionic case, and $p =
  7$ in the octonionic case. These all correspond to complex
  structures on $\mc{D}$. Then for $X,Y \in \Gamma(\mc{D})$, $\K(X,Y)$ can be
  written as: \be \K(X,Y) = \sum_{i=1}^p \gamma(I_i (X),Y) I_i.  \ee By
  extension, define $I_0$ to be the identity transformation of
  $\mc{D}$.  Now pick local orthonormal sections $\{Y_1, \cdots Y_q\}$ of
  $\mc{D}$, chosen so that $I_i Y_j$ is orthogonal to all $I_k Y_l$
  whenever $i \neq k$ or $j \neq k$. This is possible, as $\gamma$ must be
  hermitian with respect to these complex structures. Here, $q = m-1$
  for quaternionic structures, and $q =1$ for octonionic structures.

  Again, we may rewrite $\K$ as: \be \K = \sum_{i=0,j,k=1}^{i,j=p,k=q}
  -(I_j I_i Y_k)^*\otimes (I_iY_k)^* \otimes I_j.  \ee If we raise and lower all
  indexes with $\eta$ and $\gamma$, we get $\K^*$, which is \be \K =
  \sum_{i=0,j,k=1}^{i,j=p,k=q} -(I_j I_i Y_k)\otimes (I_iY_k) \otimes (I_j)^*.  \ee
  Now $\alpha \K$ involves taking the trace of $\K$ and $\K^*$ over one of
  the $\mc{D}$ components and over the $TX/\mc{D}$ components. The
  trace over the $TX/\mc{D}$ component is trivial; and if $\llcorner$
  denotes contraction between a space and its dual, \be
  \alpha \K &=& \sum_{k,o = 1}^{q} \sum_{j,r = 1}^{p} \sum_{i,l =0}^p \left( (I_j I_i Y_k)^* \llcorner (I_r I_l Y_o) \right) \left(I_j \llcorner I_r^* \right) (I_i Y_k)^* \otimes (I_l Y_o) \\
  &=& \sum_{k,o = 1}^{q} \sum_{j,r = 1}^{p} \sum_{i,l =0}^p \left( (I_j I_i Y_k)^* \llcorner (I_r I_l Y_o) \right) \delta_{jr} (I_i Y_k)^* \otimes (I_l Y_o) \\
  &=& \sum_{k,o = 1}^{q} \sum_{j = 1}^{p} \sum_{i,l =0}^p \left( (I_j I_i Y_k)^* \llcorner (I_j I_l Y_o) \right) (I_i Y_k)^* \otimes (I_l Y_o) \\
  &=& \sum_{k,o = 1}^{q} \sum_{j = 1}^{p} \sum_{i,l =0}^p \delta_{il} \delta_{ko} (I_i Y_k)^* \otimes (I_l Y_o) \\
  &=& \sum_{j=1}^p \sum_{i=0}^p \sum_{k=0}^q (I_i Y_k)^* \otimes (I_i Y_k) \\
  &=& p Id_{\mc{D}}.  \ee The $-2\beta \K$ term is the contraction of $\K$
  and $\K^*$ over both their $\mc{D}$ components; it is \be
 -2 \beta \K &=& \sum_{k,o = 1}^{q} \sum_{j,r = 1}^{p} \sum_{i,l =0}^p \left( (I_j I_i Y_k)^* \llcorner (I_r I_l Y_o) \right) \left( (I_i Y_k)^* \llcorner (I_l Y_o) \right) I_j \otimes I_r^* \\
  &=& \sum_{k,o = 1}^{q} \sum_{j,r = 1}^{p} \sum_{i,l =0}^p \left( (I_j I_i Y_k)^* \llcorner (I_r I_l Y_o) \right) (\delta_{ko} \delta_{il}) I_j \otimes I_r^* \\
  &=& \sum_{k = 1}^{q} \sum_{i =0}^p \sum_{j,r = 1}^{p} \left( (I_j I_i Y_k)^* \llcorner (I_r I_i Y_k) \right) I_j \otimes I_r^* \\
  &=& \sum_{k = 1}^{q} \sum_{i =0}^p \sum_{j,r = 1}^{p} \delta_{jr} I_j \otimes I_r^* \\
  &=& \sum_{k = 1}^{q} \sum_{i =0}^p \sum_{j = 1}^{p} I_j \otimes I_j^* \\
  &=& q (p+1) Id_{TX/\mc{D}}.  \ee Then substituting in the values for
  $p$ and $q$ gives the result.
\end{proof}

Now if $\mc{N}$ is a section of $\land^2 \mc{D}^* \otimes
(TX/\mc{D})$, we may use it in equations (\ref{al:be:k}) instead of
$\K$; in that case, define
 \be
(\alpha_{\mc{N}} \mc{M})^r_q &=& (\gamma^{jr} \gamma^{ip} \eta_{ko}) (\mc{M}_{ij}^k {\mc{N}}_{pq}^o), \\
(\beta_{\mc{N}} \mc{M})_r^k &=& - \frac{1}{2} (\eta_{or} \gamma^{ip}
\gamma^{jq}) (\mc{M}_{ij}^k {\mc{N}}_{pq}^o).  \ee
Similarly, though $\mc{K}$ defines the conformal class of $(\eta,
\gamma)$, (through the reduction to structure group $G_0 \subset G'$),
there is no reason to require that $\mc{K}$ be the Levi-bracket of the
distribution $\mc{D}$. Given $(\eta,\gamma)$ on a general manifold
with distribution $\mc{D}$ of correct dimension and co-dimension, they
define a (local) class of compatible brackets $\mc{K}$ of
quaternionic-contact or octonionic-contact type. Then the equations
(\ref{qua:con}) and (\ref{oct:con}) can be rewritten as saying that we
must find $(\eta,\gamma)$ such that for any $\mc{K}$ compatible with
them, \beqa
\label{alph:eq}\alpha_{\mc{K}} \K &=& \alpha_{\eL} \eL, \\
\beta_{\mc{K}} \K &=& \beta_{\eL} \eL, \eeqa or, more compactly, \beqa
\label{al:be:eq}
\partial^*_{\K} \K = \partial^*_{\eL} \eL.
\eeqa
It is easy to see that these equations are conformally invariant.

\begin{rem}
  It is useful to compare these equations with those defining a
  `Damek-Ricci' space (this is a subclass of Heber's metrics, see
  \cite{DamekRicci}). For any section $Z$ of $\mc{D}$, we may define
  an endomorphism $J_Z$ of $\mc{D}$ by \be \gamma (J_Z X, Y) = \eta^2
  (Z, \mc{L}(X,Y)), \ee for sections $X$ and $Y$ of $\mc{D}$. Then $X$
  is asymptotically Damek-Ricci if $J_Z^2 = -1 \eta^2(Z,Z)$. Now if
  $\{Z_j\}$ is a local orthonormal frame for $TX / \mc{D}$, then we
  may rewrite $\alpha_{\eL} \eL$ once more as \be
  \alpha_{\eL} \eL &=& \Tr_{\gamma} \big( \sum_{jk} \eta^2(Z_j,Z_k) \eta^2(Z_j, \eL) \otimes \eta^2(Z_k, \eL) \big) \\
  &=& \Tr_{\gamma} \big( \sum_{j} J_{Z_j} \otimes J_{Z_j} \big) \\
  &=& - \sum_j J_{Z_j}^2.  \ee Since $Z_j$ is normal, $J_{Z_j}^2 = -
  Id_{\mc{D}}$, and Damek-Ricci spaces must solve equation
  (\ref{alph:eq}). Similarly, for $Z$ and $Z'$ sections of $TX /
  \mc{D}$ \be -2 \eta^2 ( (\beta_{\eL} \eL) (Z), Z') &=& \Tr_{\gamma}
  \ \Tr_{\gamma} \ J_{Z'} \otimes J_{Z}.  \ee Since this must be
  symmetric, it values are determined by taking $Z = Z'$; in which
  case it is $4(m-1) \eta^2(Z,Z)$ in the quaternionic case, and
  $8\eta^2(Z,Z)$ in the octonionic one. Consequently Damek-Ricci
  spaces are special solutions of equation (\ref{al:be:eq}), as are
  any spaces that are asymptotically Damek-Ricci (i.e. spaces with
  $\eL$, $\gamma$ and $\eta^2$ such that the relation $J_Z^2 = -1
  \eta^2(Z,Z)$ holds).
\end{rem}

% An immediate consequence is that for the homogeneous models of
% quaternionic-contact and octonionic-contact manifolds, the space
% $(N,g)$ is always Einstein.

It is still somewhat unsatisfactory that there is a large class of $\mc{K}$ compatible with a given $\mc{L}$. It would be better
to have a procedure that fixes $\mc{K}$ uniquely (and hence the
quaternionic/octonionic structure, as well as $(\eta,\gamma)$).

In the quaternionic case, the dimension of $L$ (the full graded
automorphism group) is $(4m-4)^2 + 3^2 = 16m^2 - 32m + 25$, while the
group $G'$ is of dimension $(4m-4)(4m-5)/2 + 3 + 1 = 8m^2 -18m + 14$
and $G_0$ is of dimension $ (2m-2)(2m-1)/2 + 3 + 1 = 2m^2 - 3m +1$.

Looking at equation (\ref{al:be:eq}), one can see that $\alpha_{\eL} \eL$
takes values in the $\gamma$-symmetric component of $\mc{D} \otimes \mc{D}^*$,
while $\beta_{\eL} \eL$ takes values in the $\eta$ symmetric component of
$(TX/\mc{D}) \otimes (TX/\mc{D})^*$. These values are not completely
independent, however: the $\gamma$ trace of $\alpha_{\eL} \eL$ is the complete
trace of $\eL$ with itself, as is the $\eta$ trace of $-2 \beta_{\eL}
\eL$. Hence there is one extra relation, giving a total of
$(4m-4)(4m-3)/2 + 6 -1 = 8m^2 - 14 m +11$ independent equations --
just the right amount to reduce the structure group from $L$ to $G'$.

Now let us consider a slight deformation of a quaternionic-contact
structure; where $\eL = \K + \epsilon \mc{M}$. Re-writing equation
(\ref{al:be:eq}): \be
0 &=& \partial^*_{\eL} \eL - \partial^*_{\K} \K\\
&=& \epsilon \left(\partial_{\mc{M}}^* \K + \partial_{\K}^* \mc{M}\right) + O(\epsilon^2).  \ee
The $\epsilon$ term is the symmetric part of $\partial_{\K}^* \mc{M}$; so, to first
order, the requirement is that $\partial^*_{\K} \eL$ be completely
anti-symmetric. A method for fixing $\K$ is suggested by the following
lemma:

\begin{lemm}
  The equation $\partial^*_{\K} \mc{M} = 0$ consists of $16 m^2 - 32m + 24 $
  independent equations, which is exactly enough to restrict the
  structure group from $L$ to $G_0$.
\end{lemm}
\begin{proof}
  The operator $\partial^*$ takes values in $\mc{D} \otimes \mc{D}^* \oplus (TX/\mc{D})
  \otimes (TX/\mc{D})^*$. This bundle may be identified with $E(\mf{l})$,
  where $\mf{l}$ is the Lie algebra of $L$. The bracket $\K$
  defines a reduction to the structure group $G_0$ and hence
  $E_0$, a $G_0$-principal bundle. This defines the vector bundle
  $E_0(\mf{g}_0)$, with $\mf{g}_0$ the Lie algebra of
  $G_0$. The inclusion $E_0 \subset E$ defines an inclusion of this
  bundle into $E(\mf{l})$. Then paper \cite{capslo} implies that the
  image of $\partial^*$ is transverse to $E_0(\mf{g}_0)$, giving us our
  dimensionality result.
\end{proof}
So the natural candidate for fixing $\K$ would be one whose derivative
close to a quaternionic-contact structure is one where the
anti-symmetric part of $\partial_{\K}^* \mc{M}$ vanishes.

The simplest such condition is to simply require that the anti-symmetric part of $\partial_{\K}^* \mc{L}$ vanishes. Thus:
\begin{defi}[Compatibility]
  The algebraic bracket $\mc{K}$, a section of $\land^2 \mc{D}^* \otimes
  (TX/\mc{D})$, is compatible with the Levi bracket $\eL$ if:
\begin{enumerate}
\item $\mc{K}$ is of quaternionic-contact or octonionic-contact type
  -- hence the dimension and co-dimension of $\mc{D}$ is correct, and
  $\K$ defines a pair of metric $(\eta,\gamma)$ up to conformal
  transformations,
\item $\partial^*_{\K} \K = \partial^*_{\eL} \eL$,
\item $\partial^*_{\K} \eL$ is symmetric.
\end{enumerate}
\end{defi}
Now, this definition is similar, but not identical, with the condition
for non-regular two-graded geometries laid out in \cite{metwograd};
indeed, the condition there (that $\partial^*_{\K} \eL = 0$) is precisely the
infinitesimal version of the above.

Recall that in general, $\mc{L}$ is defined by an $L$-equivariant map
$f_{\eL}$ from $E$ to $W = \land^2 V_1^* \otimes V_2$. This allows us to phrase
our general result:
\begin{theo} \label{infi:ein:theo} Let $X$ be a manifold with
  distribution $\mc{D} \subset TX$ of the right dimension and
  co-dimension. Then there is an open set $U \subset W$, containing the
  standard bracket $\kappa$, such that for all points $x \in M$ where
  $f(E_x)$ intersects $U$, there exists a locally unique, continuously
  defined, choice of compatible algebraic bracket $\mc{K}_x$.
\end{theo}
\begin{proof}
  We need to prove is the existence and uniqueness of compatible
  $\mc{K}$. This is a purely algebraic construction, so we may work at
  a point. If we choose the natural bracket $\kappa$ to be fixed in $\land^2
  V_1^* \otimes V_2$, define $\theta$ as the map $\land^2 V_1^* \otimes V_2 \to \mf{s}$,
  \beqa \label{theta:equations} \theta = \frac{1}{2} \left(\partial^*_{\ell} \ell +
    \partial^*_{\kappa} \ell - \partial^*_\ell {\kappa} \right) - \partial_{\kappa}^* {\kappa}.  \eeqa Note that the
  first $\ell$ term must by symmetric, while the other two $\ell$ terms
  together are anti-symmetric, so there is no overlap between
  them. Another important fact is that the Lie algebra $\mf{g}_0$ has
  one symmetric part (the grading element) and the rest is
  anti-symmetric. We already know that the image of $\partial^*_{\ell} \ell$ is of
  co-dimension one in the symmetric part of $\mf{s}$; its image it
  precisely the part transverse to the grading element $(2Id,
  Id)$. Now $\partial^*_{\kappa} \ell - \partial^*_\ell {\kappa}$ is simply the anti-symmetric part
  of $2\partial^*_{\kappa} \ell$. We know that $2\partial^*_{\kappa} \ell$ must be transverse to
  $\mf{g}_0$, and hence so is its anti-symmetric part. Consequently
  $\theta$ maps into $\mf{l} / \mf{g}_0$.

  Now if $f_{\eL}(E_x)$ intersects the zero set of $\theta$, then there is
  a point $p \in E_x$ such that $\theta(f_{\eL}(p)) = 0$. Then if we
  define $\mc{K}_x$ by the property that $f_{\mc{K}}(p) = {\kappa}$, we
  will get the vanishing of the bundle version of equation
  (\ref{theta:equations}). Hence this $\mc{K}$ will be compatible.

  So what we need to show is that the $L$-orbit of the zero set of $\theta$
  contains an open set $U$ around ${\kappa}$. Now consider the map $\Theta: L \times
  W \to \mf{s}$, \be \Theta(s,\ell) = \theta(s \cdot \ell), \ee where $s \cdot l$ denotes the
  action of $s \in L$ on $\ell$. We wish to calculate the derivative of
  this map in the $L$ directions around the point $({\kappa},Id)$. Let $s \in
  \mf{l}$; then a little bit of calculations demonstrate that this
  derivative is \be D_{\Theta}(s)({\kappa},Id) = \partial^*_{\kappa} ( \partial_{\kappa} s) \ee where
  \be (\partial_{\kappa} s )(x,y) = {\kappa}(s(x),y) + {\kappa}(x,s(y)) - s({\kappa}(x,y)) \ee
  (see \cite{metwograd} for more details of how this is
  derived). Paper \cite{capslo} then demonstrates that $\partial_{\kappa}^* \partial_{\kappa}$
  is an invertible map from the image of $\partial^*$ to itself, with kernel
  equal to $\mf{g}_0$. An extra subtlety is needed to demonstrate that
  result, namely the vanishing of the first cohomology groups
  $H^{(1)}(\mf{g}^+,\mf{g})$ in homogeneity zero, see \cite{capslo} and
  \cite{metwograd}. But Kostant's proof of the Bott-Borel-Weil theorem
  (\cite{Kostant}) show that this is indeed the case in our situation.

  Hence, under the action of $L$, $\theta(s \cdot {\kappa})$ must trace out an open
  neighbourhood of zero in $\mf{l} / \mf{g}_0$. This property must
  extend to points $\ell$ close to ${\kappa}$ by the implicit function
  theorem, defining our set $U$.

  Now let $\ell$ be in $U$ intersected with the zero set of $\theta$. If $\ell$
  is close enough to ${\kappa}$ (possibly restricting $U$ to a smaller open
  subset), we know that if $B_\ell\subset S$ is defined such that $\theta(b \cdot \ell ) =
  0$ for all $b \in B_\ell$, then $B$ must be of same dimension as
  $\mf{g}_0$ (at least around the identity in $L$). However, if $g \in
  G_0$, then \be
  \theta(g \cdot \ell) &=& \frac{1}{2} \left(\partial^*_{g \cdot l} g \cdot \ell + \partial^*_{\kappa} (g \cdot \ell) - \partial^*_{(g \cdot \ell)} {\kappa} \right) - \partial_{\kappa}^* {\kappa} \\
  &=& \frac{1}{2} \left(\partial^*_{g \cdot \ell} g \cdot \ell + g \cdot (\partial^*_{\kappa} \ell - \partial^*_\ell {\kappa}) \right) - \partial_{\kappa}^* {\kappa} \\
  &=& \frac{1}{2} (\partial^*_{\ell} \ell) - \partial_{\kappa}^* {\kappa} = 0,\\
  \ee since $g$ is a conformal transformation, commutes with $\partial^*$,
  and $\partial^*_{\kappa} \ell - \partial^*_\ell {\kappa} = 0$ by the assumption $\theta(\ell) = 0$.

  Hence around the identity, dimension count implies that $B_\ell$ is
  precisely the group $G_0$. Action by $G_0$ preserves ${\kappa}$, so does
  not affect the value of $\mc{K}_x$. Consequently, the choice of
  $\mc{K}_x$ is locally unique for $\ell \in U$.

\end{proof}

\begin{rem}
  As noted before, the condition that $\partial^*_{\K} \eL$ be symmetric can
  be replaced with any other condition that approximates the one above
  to first order. There are more natural candidates for that --
  involving, for instance, the decomposition of the partial trace
  $\eL_{ij}^k \eL_{lo}^r \gamma^{il}$ into irreducible $G_0$ components,
  and the vanishing of one of these components. But since we're been
  unable to find a direct use of such a result (it affect the
  curvature of the asymptotically Einstein metric, but it's not clear
  exactly how), we've stuck with the simpler condition in this paper.
\end{rem}

\section{Construction of the Einstein metrics}
\label{sec:constr-einst-metr}

In this section, we prove theorem \ref{th:main}, along the lines of
\cite{Biq00}. Because we restrint to the case of small deformations of
the model hyperbolic metric, we are able to give a short direct proof,
in which the main step is a uniform estimate for the norm of the
inverse of the linearization.

We start with the quaternionic or octonionic hyperbolic space $M$, whose
metric in polar coordinates is expressed in both cases by
\begin{equation}
  \label{eq:1}
  g_0 = dr^2 + \sinh^2(\tfrac r2)\gamma_0 + \sinh^2(r) \eta_0^2 .
\end{equation}
Here $\eta_0$ is a 1-form on $S^{4m-1}$ (resp. $S^{15}$) with values in
$\setR^3$ (resp. $\setR^7$), and $\gamma_0$ is the induced metric on the
$4(m-1)$-dimensional (resp. $8$-dimensional) distribution $\cD_0$ of
$S^n$.

We will need the mean curvature $H_0(r)=\partial_r\log v$ of the spheres
$r=\mathrm{cst}$, where $v$ is the volume element. It is given by
$H_0(r)=2(m-1)\coth(\frac r2)+3\coth(r)$ in the quaternionic case, or
$H_0(r)=4\coth(\frac r2)+7\coth(r)$ in the octonionic case. Also we
note $$\cH=\lim_{r\to\infty}H_0(r) $$ the limit at infinity, so that
$\cH=2m+1$ in the quaternionic case and $\cH=11$ in the octonionic
case.

Suppose that we have now a small perturbation $\cD$ of the
distribution $\cD_0$. From proposition \ref{prop:constr-ae} we have
constructed $(\gamma,\eta)$ with $\cD=\ker \eta$ such that the metric
\begin{equation}
  \label{eq:2}
  g_\cD = dr^2 + \sinh^2(\tfrac r2)\gamma + \sinh^2(r) \eta^2
\end{equation}
is asymptotically Einstein :
\begin{equation}
  \label{eq:3}
  \Ric(g_\cD)-\lambda g_\cD = O(e^{-\frac r2}) ,
\end{equation}
with $\lambda=-m-2$ (resp. $\lambda=-9$). Here the norms are with respect to
$g_\cD$. Actually, in the proof of proposition \ref{prop:constr-ae},
we proved more, that is there is a developpement for the curvature,
\begin{equation}
  \label{eq:18}
  R = R_0 + e^{-\frac r2} R_1 + e^{-r} R_2 + \cdots ,
\end{equation}
where the terms $R_i$ do not depend on $r$, the term $R_1$ depends on
one derivative of the bracket $\mc{L}$ on the boundary, and the other
terms depend on two derivatives of $\mc{L}$. This immediately implies
\begin{equation}
  \label{eq:4}
  |\nabla^k(\Ric(g_\cD)-\lambda g_\cD)| \leq c_k e^{-\frac r2}  \text{ for all }k,
\end{equation}
where $c_k$ can be made small if $\mc{L}$ is $C^{k+2}$ close to the
standard bracket.

Of course, the formula (\ref{eq:2}) does not give a smooth metric at
the origin. To remedy this, we choose a cutoff function $\chi(r)$, such
that $\chi(r)=1$ for $r\geq R+1$ and $\chi(r)=0$ or $r\leq R-1$. Then we define
\begin{equation}
  \label{eq:5}
  g = \chi g_\cD + (1-\chi) g_0 .
\end{equation}
The metric $g$ is a global filling of $(\gamma,\eta)$ in the ball.

The first observation is that the metrics $g_\cD$ have uniform geometry:
\begin{lemm}\label{lem:uniform-geom}
  Suppose $k\geq 2$. For $\cD$ varying in a fixed $C^{k+1}$ neighbourhood
  of $\cD_0$, the sectional curvature of $g$ is negative, the
  curvatures of $g$ and their $(k-2)$ covariant derivatives are
  uniformly bounded.
\end{lemm}
\begin{proof}
  A $C^{k+1}$ control of $\cD$ gives a $C^k$ control of the conformal
  metric $(\eta,\gamma)$, since one derivative is needed to calculate the Levi
  bracket and $(\eta,\gamma)$ is then obtained as the solution of algebraic
  equations. Therefore we have a $C^k$ control on the coefficients of
  $g$. The lemma then follows from the form (\ref{eq:18}) of the
  curvature.
\end{proof}

This implies that balls for the metrics $g_\cD$ are uniformly
comparable with Euclidean balls. Then the Hölder norm of a function
$f$ is defined as the supremum of the Hölder norms of $f$ on each ball
of radius $1$.

The analysis of the Einstein equation requires the use of weighted
Hölder spaces. Our weight function will be
\begin{equation}
  \label{eq:6}
  w(r) = \cosh(r)^\delta
\end{equation}
and we then define the weighted Hölder space $C^{k,\alpha}_\delta=w^{-\delta}C^{k,\alpha}$.
Of course, from the initial estimate (\ref{eq:3}), the weight we are
interested in is $\delta=\frac 12$.

As is \cite[chapter I]{Biq00}, the Einstein metric will be constructed as a
solution $h$ of the equation
\begin{equation}
  \label{eq:7}
  \Phi^g(h):=\Ric(h)-\lambda h+\delta_h^*(\delta_gh+\tfrac 12 d \Tr_gh) = 0,
\end{equation}
and we require that $h$ is asymptotic to $g$ in the sense that
\begin{equation}
  \label{eq:8}
  h-g \in C^{2,\alpha}_{1/2}.
\end{equation}
Indeed, by \cite[lemma I.1.4]{Biq00}, a solution $h$ of $\Phi^g(h)=0$
then satisfies $\delta_gh+\tfrac 12 d \Tr_gh=0$ and $\Ric(h)=\lambda h$. Given
lemma \ref{lem:uniform-geom} (in particular, the negative curvature of
$g$ implies that the linarization of $\Phi^g$ has no $L^2$ kernel), the
proof in \cite{Biq00} applies and proves that if the data $(\gamma,\eta)$ is
sufficiently close to $(\gamma_0,\eta_0)$ in $C^{2,\alpha}$ norm, that is if $\cD$
is sufficiently close to $\cD_0$ in $C^{3,\alpha}$ norm, then one can find
a solution $h$ of (\ref{eq:7}), if one has a uniform bound on the
inverse of the linearization of $\Phi^g$. This is provided by:

\begin{lemm}\label{lemm:P_invertible}
  Suppose that $\frac 12(\cH-\sqrt{\cH^2-8})<\delta<\frac
  12(\cH+\sqrt{\cH^2-8})$. For $\cD$ sufficiently close to $\cD_0$ in
  $C^{3,\alpha}$ norm, the linearization
  $P_g=d_g\Phi^g:C^{2,\alpha}_\delta(\Sym^2T^*M)\to C^\alpha_\delta(\Sym^2T^*M)$ is invertible
  and the norm of the inverse is uniformly bounded.
\end{lemm}

From the value of $\cH$, we check that the weight $\delta=\frac 12$ indeed
satisfies the hypothesis, so theorem \ref{th:main} follows from the
lemma.

So we now concentrate on the proof of the lemma. One has
$$P_g = \tfrac 12 \nabla^*\nabla - \Rr_g . $$
The property of the curvature term $\Rr_g$ we need is the following
\cite[lemmas I.4.1 and I.4.2]{Biq00}: for the hyperbolic metric $g_0$,
the largest eigenvalue of $\Rr_{g_0}$ is equal to $1$ (instead of $4$
in \cite{Biq00}, because here we normalize here the sectional
curvature of $g_0$ in $[-1,-\frac 14]$ instead of $[-4,-1]$). This
immediately implies that, for $\cD$ close enough to $\cD_0$ in $C^3$
norm, one has
\begin{equation}
  \label{eq:9}
  \Rr_g \leq 1 + \epsilon .
\end{equation}

For the function $w$ depending on $r$ only, one has
\begin{equation}
  \label{eq:10}
  \Delta w = - \partial_r^2w - H(r) \partial_rw, 
\end{equation}
where $H(r)=\partial_r\log v$ is the mean curvature. For the metric $g$ given by
(\ref{eq:5}), the mean curvature $H(r)$ coincides with $H_0(r)$ for
$r\geq R+1$ or $r\leq R-1$, and for $R-1\leq r\leq R+1$ we get $|H(r)-H_0(r)|\leq \epsilon$
if we suppose $(\gamma,\eta)$ close enough to $(\gamma_0,\eta_0)$.

An easy calculation gives, for the hyperbolic metric,
\begin{equation}
  \label{eq:11}
  -\frac{\Delta w}w-2\frac{|dw|^2}{w^2} =
     \delta \left( \cH-\delta + \frac{\dim \cD}{2\cosh r} + \frac{\delta+1}{\cosh^2r} \right).
\end{equation}
It follows that, for the metric $g$, if $\cD$ is sufficiently close to
$\cD_0$,
\begin{equation}
  \label{eq:12}
  -\frac{\Delta w}w-2\frac{|dw|^2}{w^2} \geq \delta (\cH-\delta-\epsilon).
\end{equation}

Using this property of the weight function $w$, we can now establish
lemma \ref{lemm:P_invertible} using the maximum principle. From Kato's
inequality, 
$$ \langle u,\nabla^*\nabla u\rangle = |u|\Delta|u| + |\nabla u|^2 - |d|u| |^2 \geq |u|\Delta|u| . $$
Using the formula
\begin{align*}
 w\Delta|u| &= \Delta(w|u| )-w|u|\big(\frac{\Delta w}w+2\frac{|dw|^2}{|w|^2}\big)
          +2\big\langle\frac{dw}w,d(w|u| )\big\rangle\\
       &\geq \Delta(w|u| )+\delta (\cH-\delta-\epsilon)w|u|+2\big\langle\frac{dw}w,d(w|u| )\big\rangle
\end{align*}
it follows from (\ref{eq:9}) that
\begin{equation}
  \label{eq:13}
  w|Pu| \geq \tfrac 12 \Delta(w|u| )+ \big( \tfrac 12 \delta(\cH-\delta-\epsilon) - 1 - \epsilon \big) w|u| + \big\langle\frac{dw}w,d(w|u| )\big\rangle .
\end{equation}
Let $A=\tfrac 12 \delta(\cH-\delta-\epsilon) - 1 - \epsilon$. If $\delta$ satisfies the hypothesis
of lemma \ref{lemm:P_invertible}, then one can choose $\epsilon$ sufficiently
small so that $A>0$. Then by the maximum principle applied to $w|u|$,
it follows that
\begin{equation}
  \label{eq:13b}
  \sup (w|u| ) \leq A^{-1} \sup (w|Pu| ).
\end{equation}
(A priori we cannot apply the maximum principle to $w|u|$ since it has
not to go to zero at infinity, but we can apply it for $w=(\cosh
r)^{\delta'}$ for any $\delta'<\delta$; then taking $\delta' \to \delta$ gives the estimate).

From this estimate, it is immediate that if $v\in C^\alpha_\delta$, then one can
solve $Pu=v$ with $u\in C^0_\delta$ and $\|u\|_{C^0_\delta}\leq A^{-1}\|v\|_{C^0_\delta}$. It
remains to obtain a bound on higher derivatives, but from the uniform
geometry lemma \ref{lem:uniform-geom}, applying the usual elliptic
estimate in each ball, one obtains a constant $C$ such that
$$ \|u\|_{C^{2,\alpha}_\delta} \leq C \big( \|Pu\|_{C^\alpha_\delta} + \|u\|_{C^0_\delta} \big) \leq
C(1+A^{-1})\|Pu\|_{C^\alpha_\delta} $$
which is the required estimate.

\begin{rem}
  The previous lemma does not give an optimal interval of weights for
  the isomorphism. In \cite{Biq00} the optimal interval for $g_0$ is
  calculated; using microlocal analysis, it is proved in \cite{BM}
  that the same interval holds if the distribution $\cD$ is
  quaternionic-contact (the regular case). In general, the optimal
  interval may depend on the supremum of the eigenvalues of the
  curvatures $\Rr_x$, where $\Rr_x$ is the curvature of the homogeneous
  Einstein model attached to the point $x$ of the boundary.
\end{rem}

\bibliographystyle{amsalpha}
\bibliography{quat}

\end{document}